\documentstyle{amsppt}                                                      

\baselineskip=24pt
\pagewidth{125mm}
\pageheight{175mm}
\NoRunningHeads
\loadeusm

\magnification = 1200

\font\srp=eusm10

\define\BP{\bigpagebreak}

\define\q{\qquad}
\define\qq{\qquad\qquad}

\define\SP{\smallpagebreak}

\define\a{\,$}
\define\A{\;$}

\define\bB{\bold B}
\define\bC{\bold C}

\define\CF{\Cal F}
\define\CL{\Cal L}

\define\CH{\Cal H}

\define\CS{\Cal S}

\define\Con{C_0^\infty}

\define\CX{\Cal X}

\define\f{$\,}

\define\gr{\nabla}

\define\M{\Vert}

\define\Mtt1{{\Tilde {\Tilde M}}_1}

\define\pa{\partial}

\define\pd#1#2{\dfrac{\partial#1}{\partial#2}}

\define\Ro{\bold R}
\define\R3{\bold R^3}
\define\Rn{\bold R^n}

\define\SF{\text{\srp F}}

\define\T{\tag}

  
\topmatter
\title Resolvent Estimates of the  Dirac Operator \endtitle
\author Chris Pladdy$^1$, Yoshimi Sait\B o$^2$, and Tomio Umeda$^3$
                                                               \endauthor
\affil  $^{1,2}$ Department of Mathematics \\ 
          University of Alabama at Birmingham \\ 
          Birmingham, Alabama 35294 \\
          U. S. A. \\
                             and \\ 
        $^3$ Department of Mathematics \\
          Himeji Institute of Technology \\
          Himeji 671-22 \\ Japan                        \endaffil
\endtopmatter

\newpage

\document

\heading      {\bf \S1. Introduction}              \endheading 

\SP

        The present paper is concerned with the Dirac operator
$$
  H = -i \sum_{j=1}^3 \alpha_j \frac{\pa}{\pa x_j} + \beta + Q(x),  \T 1.1
$$
     where \f i = \sqrt{-1} $, \f x = ( x_1, x_2, x_3 ) \in \R3 \A   and
     \f \alpha_j $, \f \beta \A are
     the Dirac matrices, i.e.,  $4 \times 4$ 
     Hermitian matrices satisfying the anticommutation relation
$$
 \qq
  \alpha_j\alpha_k + \alpha_k\alpha_j = 2\delta_{jk} \, I 
                     \qq (j, k = 1, 2, 3, 4)                         \T 1.2
$$
     with the convention $\,{\alpha}_4 = \beta\;$ , \ $\, \delta_{jk}  \;$ 
     being  Kronecker's delta and $\,I\;$  being the $4 \times 4$ identity 
     matrix. The potential  $\, Q(x) \;$ is a  $4 \times 4$ 
     Hermitian matrix-valued function, which is usually
     assumed  to diminish at infinity.                               \par

        The limiting absorption principle for the operator $\, H \; $
     was first discussed by Yamada\;[13]. 
     As a result, the existence of the extended resolvents 
     $\,R^{\pm} (\lambda)\;$ was assured (see Theorem 2.2 in section 2 
     below). The extended resolvents   $\,R^{\pm} (\lambda)\;$
     play important roles in spectral and scattering theory for the 
     operator $\, H \;$ (see [13] and [14]).                           \par

        The aim of this paper is to investigate the asymptotic behavior of
     $\,R^{\pm} (\lambda)\;$ as $\, |\lambda|  \to  \infty \,$.
     As Yamada\; [15] pointed out, the operator norm in 
     $\, \bB(\CL_{2,s}, \CL_{2,-s}) \;$ of the extended resolvents   
     $\,R_0^{\pm} (\lambda)\;$ of the free Dirac operator
$$
        H_0 = -i \sum_{j=1}^3 \alpha_j \frac{\pa}{\pa x_j} + \beta
                                                                   \T 1.3
$$
     cannot approach zero as $\, |\lambda|  \to  \infty\,$. 
     (The definition of $\, \bB(\CL_{2,s}, \CL_{2,-s}) \;$  is found
     below in the introduction.)  This means that boundedness of the 
     operator norm of  $\,R_0^{\pm} (\lambda)\;$ is possibly the
     best that one can show. Indeed, one of our main  results is that 
     the operator norm of $\,R_0^{\pm} (\lambda)\;$ stay  bounded as 
     $\, |\lambda|  \to  \infty\,$ (see Theorem 2.4 below). However, we
     also show that  $\,R_0^{\pm} (\lambda)\;$
     converge strongly  to  $0$ as $\, |\lambda|  \to  \infty\,$ (see 
     Theorem 2.5). Our results indicate that   the
     extended resolvents of Dirac operators decay much more slowly
     than those of Schr\"odinger operators.
     (Compare Theorems 2.4 -- 2.6 with Theorem 4.1.)                    \par

        We now introduce the notation which will be used in this paper. 
     For $\,x \in \R3\,$, $\,|x|\,$ denotes the Euclidean norm of $\,x\,$ 
     and
$$
     \langle x \rangle = \sqrt{1 + |x|^2}.                          \T 1.4
$$
     For \f s \in \Ro $, we define the weighted Hilbert spaces 
     \f L_{2,s}(\R3) \a and \f H^1_s(\R3) \a by           
$$
    L_{2,s}(\R3) = \{ f \ / \ {   \langle x \rangle }^s f
                                         \in L_2(\R3) \},          \T 1.5     
$$
     and
$$
    H^1_s(\R3) = \{ f \ / \
      {   \langle x \rangle }^s
        \pa_x^{\alpha}f \in L_2(\R3), |\alpha|
                                           \le 1 \},           \T 1.6
$$
     where $\,\alpha = (\alpha_1, \alpha_2, \alpha_3)\,$ is a multi-index,
     $\,|\alpha| = \alpha_1 + \alpha_2 + \alpha_3$, and 
$$
     \big( \frac{\pa}{\pa x} { \big) }^{\alpha}
       = \frac{\pa^{|\alpha|}}
             {\pa x_1^{\alpha_1}\pa x_2^{\alpha_2}\pa x_3^{\alpha_3}}.
                                                                \T 1.7
$$ 
     The inner products and norms in $\,L_{2,s}(\R3)\,$ and $\,H^1_s(\R3)\,$ 
     are given by
$$
\left\{ \aligned
       &(f, g)_s = \int_{\R3} { \langle x \rangle }^{2s}
              f(x)\overline{g(x)}\, dx, \\
       &\M f \M_s = \big[(f, f)_s \big]^{1/2},  
\endaligned \right.                                               \T 1.8
$$ 
     and   
$$
\left\{ \aligned
       & (f, g)_{1,s} = \int_{\R3} { \langle x \rangle }^{2s}
          \big[\,  f(x)\overline{g(x)}\, + \,
                  \gr f(x) \cdot
                         \overline{ \gr g(x)} 
                                    \,  \big]\, dx,   \\
       & \M f \M_{1,s} = \big[(f, f)_{1,s} \big]^{1/2},    
\endaligned \right.                                               \T 1.9
$$ 
     respectively. The spaces \f \CL_{2,s} \a and \f \CH^1_s \a are 
     defined by
$$
\left\{ \aligned
     & \CL_{2,s} = \big[ L_{2,s}(\R3) \big]^4, \\
     &  \CH^1_s = \big[ H^1_s(\R3) \big]^4, 
\endaligned \right.                                              \T 1.10
$$
     i.e., \f \CL_{2,s}\,$ and \f \CH^1_s\,$ are direct sums of the 
     Hilbert spaces $\,L_{2,s}(\R3) \a and \f H^1_s(\R3)$, respectively. The 
     inner products and norms in \f \CL_{2,s} \a and \f \CH^1_s \a 
     are also denoted by \f (\ , \ )_s$, \f \M \ \M_s \a and 
     \f (\ , \ )_{1,s}$, \f \M \ \M_{1,s}$, respectively. When \f s = 0$, 
     we simply write
$$
\left\{ \aligned
       & \CL_2 = \CL_{2,0}, \\
       & \CH^1 = \CH^1_0.   
\endaligned \right.                                                \T 1.11
$$
     For a pair of Hilbert spaces \f X \A and \f Y $, \f \bB(X, Y) \A 
     denotes the Banach space of all bounded linear operators from \f X \a 
     to \f Y$, equipped with the operator norm
$$
      \M T \M = \sup_{x\in X\backslash \{0\}} \M Tx \M_Y / \M x \M_X,
                                                                  \T 1.12          
$$
     where $\,\M \ \M_X\,$ and $\,\M \ \M_Y\,$ are the norms in \f X \a and 
     \f Y$.
     For \f T \in \bB(\CL_{2,s}, \CL_{2,t})$, its operator norm will be
     denoted by \f \M T \M_{(s, \, t)} $.                           \par

     We now sketch the contents of the paper. In section 2, we state
     the main theorems. For the reader's convenience, we
     reproduce Yamada's arguments\;[15] in section 3. 
     In section 4, we make a brief review of resolvent estimates
     for 
     Schr\"odinger operators which will be 
     used in the proof
     of Theorem 2.4. In section 5, we establish some boundedness
     results  for 
     pseudodifferential operators acting in the weighted Hilbert spaces,
     the results  on which
     the proof of Theorem 2.4 is based. We give the proofs of Theorems
     2.4 and 2.5 
     in sections 6 and 7 respectively.  In section 8, we give the 
     proof of  
     Theorem 2.6.                                                 \par
     
     Finally, we would like to mention that Pladdy, Sait\B o 
     and Umeda\,[6]
     is an announcement for this work.  Also, We would like to mention
     that we can establish resolvent estimates for 
     relativistic Schr\"odinger operators
     $\,  \sqrt{-\Delta + m^2} \, + \, V(x)$, the estimates which are
     similar to those of 
     the Dirac operators. Discussions about the resolvent estimates
     for the relativistic Schr\"odinger operators 
     will appear elsewhere.

     The present work was done while the last author (T.U.) was
     visiting the Department of Mathematics,
     the University of Alabama at Birmingham for the 1992--93 academic
     year. He would like to
     express his sincere gratitude to
     the members of  the department for their
     warm hospitality. He also
     would like to   thank Himeji Institute of Technology for
     allowing him to take a year's leave of absence.               \par


\vskip 24pt

\heading      {\bf \S2. Main results}              \endheading 

\SP

    We begin with the selfadjointness of the free Dirac operator 
    $\,H_0$. It is known that $\,H_0\;$ restricted on  
    $\, \big[ \Con(\R3) \big]^4\;$
    is essentially 
    selfadjoint in  $\, \Cal L_2 \;$ and its selfadjoint extension,
    which will be denoted by
    \f H_0 \A  again, has the domain \f \Cal H^1 $.             \par

    We impose the following assumption on the potential.
\BP
        
        {\bf Assumption 2.1.}        
\SP

        (i) \f Q(x) = (q_{jk}(x)) \A is a \f 4 \times 4 \a Hermitian 
     matrix-valued \f C^1 \a function on \f \R3$;
     
        (ii) There exist positive constants \f \epsilon \A  and \f K \A  such 
     that
$$
        { \langle x \rangle }^{1+\epsilon}
              |q_{jk}(x)| + 
       \sum_{\ell=1}^3 |\frac{\pa q_{jk}}{\pa x_{\ell}}(x)| \le K  \T 2.1              
$$
     for $\,j, k = 1, 2, 3, 4$.        
\BP

     Assumption 2.1 is essentially the same one as Yamada  made  
     in [13]. He needs the first derivatives of \f q_{jk} \A 
     bounded in order to assure that the Dirac operator \f H \A
     has no embedded eigenvalues in its essential spectrum;
     see [13, Proposition 2.5].

        Under Assumption 2.1 the multiplication operator $\,Q = Q(x)\times\;$
     is a bounded selfadjoint operator in $\,\CL_2$. Hence, by the 
     Kato-Rellich theorem (Kato[4], p.287), $\,H\;$ restricted on  
     $\,[\Con(\R3)]^4\,$ is also essentially selfadjoint in \f \CL_2 \A 
     and its selfadjoint extension, which will be denoted by \f H \A again, 
     has the same domain \f \CH^1 \A  as \f H_0$.  We write
$$
             R_0(z) = (H_0 - z)^{-1},                         \T 2.2         
$$
     and
$$
             R(z) = (H - z)^{-1}.                             \T 2.3           
$$
     As we mentioned in the introduction, 
     the limiting absorption principle
     holds for the Dirac operator \f H$.  We note that in Theorem 2.2
     below $\, R(z) \;$ is regarded 
     as an operator  belonging to \f \bB(\CL_{2,s}, \CL_{2,-s})$.
\BP
      
        {\bf Theorem 2.2}\ {\it (Yamada\;[13]).
\SP
        
        Suppose that Assumption 2.1 is satisfied and let $\,s > 1/2$. Then
     for  any                                             
     \break \f \lambda \in (-\infty, -1) \cup (1, \infty)$, there 
     exist the extended resolvents \f R^{\pm}(\lambda) \in 
     \bB(\CL_{2,s}, \CL_{2, -s}) $\break such that 
     for any   \f  f  \in   \CL_{2,s} $ }
$$
    R(\lambda \pm i\mu) f   \longrightarrow
                 R^{\pm}(\lambda) f  \;\;\;\;
                   \text{\it in } \; \CL_{2, -s} 
$$
   {\it  as \f \mu  \downarrow 0  $. 
      Moreover, for \f f \in \CL_{2,s}$, \f R^{\pm}(\lambda)f \A is
     an \f \CL_{2,-s}$--valued  continuous function on 
     \f (-\infty, -1) \cup (1, \infty)$.}                                                                            
\BP

   { \bf Remark 2.3. } 
\SP
     
\item{(i)}  Actually, Yamada\;[13] proved that \f R^{\pm}(\lambda) \A
     belong to \f \bB(\CL_{2,s}, {\CH}^1_{-s})$.  In particular,
     \f R^{\pm}(\lambda) \A belong to 
    $ \bB(\CL_{2,s}, \CL_{2, -s}) $,  which is suitable for our purpose.
\item{(ii)}  Note that the conclusions of Theorem 2.2 are valid,
      in particular, for the resolvent of the free Dirac operator 
      \f H_0 $.
\BP

        We now state the main theorems, which are concerned with the 
     asymptotic behavior of the extended resolvents \f R_0^{\pm}(\lambda) \A
     of the free Dirac operator \f H_0$.

\vskip 30pt

        {\bf Theorem 2.4.}
\SP
     
        \it Let $\,s > 1/2$. 
     Then 
$$
        {  \M R_0^{\pm}(\lambda) \M }_{(s, -s)}
             = O(1)  \qq (|\lambda| \to \infty).
$$
\rm
\vskip 24pt

        As we shall see later,
          $\,{\M R_0^{\pm}(\lambda) \M}_{(s, -s)}\;$
           cannot 
     be small no matter how $\,|\lambda|\,$ is large. In this sense the 
     estimate in Theorem 2.4 is best possible. However, 
       $\, R_0^{\pm}(\lambda) \;$
      do become small  
     in strong operator topology  as $\,|\lambda|\;$ gets  large.
     In fact, we have
\BP

        {\bf Theorem 2.5.}
\SP
     
        {\it Let \f s > 1/2 $. 
     Then \f R_0^{\pm}(\lambda) \A converge strongly to \f 0 \A  as 
     \f |\lambda| \to \infty$, i.e., 
         for any \f f \in \CL_{2,s}$             
$$
         R_0^{\pm}(\lambda)f \longrightarrow  0    
         \;\;\; \hbox{ in } \, \CL_{2,-s}                        \T 2.4                           
$$
    as \f |\lambda| \to  \infty$.}
\BP

  Based on Theorems 2.4 and 2.5, the Dirac operator with a small
     coupling constant can be handled; we  can use the Neumann
     series expansion. Let
$$
       H_t =
        -i \sum_{j=1}^3 \alpha_j \frac{\pa}{\pa x_j} + \beta + t \,  Q(x),
$$
     where \f t  \A  is a real number.  The extended resolvents of
    \f  H_t  \A  will be denoted by 
     \f R_t^{ \pm} (\lambda)$.
      Then we have the following 

\BP

        {\bf Theorem 2.6.}
\SP

        \it Suppose that  \f Q(x) \A  satisfies Assumption 2.1 and
      let   \f s > 1/2 $. Then there exist
     constants \f t_0 > 0 \A  and \f C > 0 \A  
     such that    
       for every \f t \A with
             \f | t  | \, \le \,  t_0$
  \item{\rm (i)}  $       \sup_{|\lambda| \ge 2}
         { \M  R_t^{\pm} (\lambda) \M }_{(s, -s)} \,
             \le  \, C $,
   \item{\rm(ii)}  \f R_t^{ \pm} (\lambda) \A  converge strongly to  
     \f 0 \A as \f | \lambda |  \to   \infty $.  \rm

%
%
 
\vskip 24pt

\heading      {\bf \S3. Yamada's counterexample}              \endheading 

\SP

     In this section, we shall  reproduce Yamada's arguments\;[15] to
     show that                                              \break
     $\,{  \M R_0^{\pm}(\lambda) \M }_{(s, -s)}\;$  cannot
     converge to  $0$   as  $| \lambda |  \to  \infty$.
     In Proposition 3.1 below and in the
     rest of the paper,  $\CS  (\R3)$  denotes
     the Schwartz space
      of rapidly decreasing functions on $ \R3 $.
      
\vskip 12pt

        {\bf Proposition 3.1 \rm ( \it  Yamada [15] ). }
\SP

        \it  There exists a sequence $\{  h_n \}{}_{n=1}^{\infty}
         \subset [\CS  (\R3) ]^4 $ 
         such that
\item{\rm(i)} \  $\sup_n  {\M h_n \M}_s  \, < \,   + \infty\;\;$
         for every $s>0$,
\item{\rm(ii)} \ $\lim_{n \to \infty} \, (R_0^{\pm}(n+2)h_n, \, h_n)_0
          \,  \not=  \, 0$.   \rm

\vskip 12pt

     It follows from Proposition 3.1 that for any $\, s > 1/2 $, 
      $\,{  \M R_0^{\pm}(\lambda) \M }_{(s, -s)}\;$  cannot
     converge to  $0$   as  $ \lambda   \to  \infty$.
     In fact, the inequality
$$
   |(R_{0}^{\pm}(n+2)h_n, h_n)_0|  \,  \le \,
       { \M R_{0}^{\pm}(n+2) \M }_{(s, -s)}  \, {\M h_n \M_s}^2,   
$$
     together with Proposition 3.1, implies that 
$$
     \liminf_{n \to \infty} { \M R_{0}^{\pm}(n+2) \M }_{(s, -s)} \,  > \, 0.
$$
%

\newpage

        {\bf Remark 3.2. }
\SP

     One can also show that for any $\, s>1/2 $, 
         $\,{  \M R_0^{\pm}(\lambda) \M }_{(s, -s)}\;$  cannot
     converge to  $0$   as  $ \lambda   \to  - \infty$.
     See Yamada\;[15].

\vskip 12pt

   Throughout (and only in) this section, we assume that
$$
 \beta = \pmatrix
          1 & 0 & 0 & 0  \\
          0 & 1 & 0 & 0  \\
          0 & 0 & -1 & 0  \\
          0 & 0 & 0 & -1
    \endpmatrix.                                                 \T 3.1
$$
     This causes no loss of generality. Indeed,  $H_0\;$ with
     any Dirac matrices
     is unitarily equivalent to 
     $\, H_0 \;$ with 
     the Dirac matrices \f {\alpha}_j \A  and \f \beta \A 
     of the form  (3.1).                                           \par

     We will give the proof of Proposition 3.1 with a series of lemmas.

\vskip 12pt

        {\bf Lemma 3.3. }
\SP

   \it     Let \f \varphi \A be a  real-valued $C^1$--function 
        defined on $[-1, \, 1]$.  Then
$$
 \lim_{\mu \downarrow 0} \, \int_{-1}^1 \,
    \frac{\varphi ( \sigma )}{ \sigma \,  \mp  \, i \mu} \, d\sigma
      \; = \;  
         \pm   i \pi \varphi (0) \, + \, 
         \int_{-1}^1 \big\{ \int_0^1   { \varphi }^{\prime}
          ( \sigma  \theta )  \,  d\theta \,  \big\} \, d\sigma. 
$$
\vskip 6pt
       Proof. 
\rm   \ It is easy to see that 
$$
 \int_{-1}^1 \,
    \frac{\varphi ( \sigma )}{ \sigma \, \mp \, i \mu} \, d\sigma
      \; = \;  
  \varphi (0) \, 
 \int_{-1}^1 \,
    \frac{1}{ \sigma \, \mp \, i \mu} \, d\sigma   \; + \;
    \int_{-1}^1 \,
    \frac{\varphi ( \sigma ) \, - \, \varphi (0 )}
    { \sigma \, \mp \, i \mu} \, d\sigma.                             \T 3.2
$$
     Noting that 
$$
\varphi ( \sigma ) \, - \, \varphi (0 ) \, = \,  
              \sigma \, \int_0^1 \, { \varphi }^{\prime}
                 (\sigma \theta) \,  d\theta,
$$
     and taking the limit of (3.2) as  $\,  \mu  \downarrow 0 $,
     we get the desired conclusion.   $ \square$
\vskip 12pt

        {\bf Lemma 3.4. }
\SP

   \it     Let \f \varphi \A be a  real-valued $C^1$--function 
        defined on $[-1, \, 1]$   and suppose that 
        $\,\varphi\;$  is an 
        even function.  Then   
$$
   \int_{-1}^1 \big\{ \int_0^1   { \varphi }^{\prime}
          ( \sigma  \theta )  \,  d\theta  \,  \big\} \,
                                          d\sigma   \, = \, 0. 
$$
\vskip 6pt
Proof.
\rm
      Since  \f \varphi \A is an even function, we see that
      \f {\varphi}^{\prime} \A  is an 
      odd function.  Then    
$\,  \int_0^1   { \varphi }^{\prime}
          ( \sigma  \theta )  \,  d\theta   \;$   is
          also an odd function of  \f  \sigma $,
     of which integral from $-1$ to $1$ is equal to $0$.   
                                                    $\square$
\vskip 12pt
    To prove Proposition 3.1,  we shall construct the sequence
    $\, \{ h_n \} \;$  in the 
    following manner: First choose
    an even function  $\, \varphi  \in  C_0^{\infty} ( \Ro ) \;$
    so that
$$
\text{supp}[ \,  \varphi ] \,  \subset \,  (-1, \, 1)              \T 3.3
$$
   and
$$
     \varphi (0) \, = \, 1.                                       \T 3.4
$$
    Next, define  \f a_n \in \CS ( \R3 ) \A by
$$
\widehat a_n ( \xi ) \, = \,
    \frac1{ | \xi |}  \,  \varphi (
        \, \langle \xi  \rangle  - n -2 \, )  \qq (n=1, 2,\cdots),  \T 3.5
$$
   where  \f \widehat a \, = \, \SF a \A 
   is the Fourier transform of \f a $:
$$
  \widehat a ( \xi ) \, = \, [ \SF a ] (\xi) \, = \, \int_{\R3}
     \, e^{-ix \cdot \xi} \, a(x) \, dx.
$$
    Later we will also 
    use the inverse Fourier transform which is given by
$$
   [{ \SF}^{\,-1} b ] (x) \, = \,  (2 \pi )^{-3} \int_{\R3}
     \, e^{ix \cdot \xi} \, b(\xi) \, d\xi.
$$
     Note that 
$$
 \text{supp} [\, \widehat a_n ]  \, \subset \,
      \big \{ \;  \xi \in \R3 \;  \big/
              \;\;  
               \sqrt{n(n+2)}  \le | \xi |  \le
                  \sqrt{n^2 + 6n + 8} \; \big\}.                  \T 3.6
$$
   Finally, define \f h_n \, \in \, [ \CS ( \R3 ) ]^4 \A  by
$$
h_n \, = \, \pmatrix
            a_n \\
            0 \\
            0 \\
            0
         \endpmatrix   
             \qq ( n=1,2, \cdots ).                                \T 3.7
$$
\vskip 12pt

        {\bf Lemma 3.5. }
\SP

   \it   For any \f s > 0 $
$$
 \sup_n \, { \M h_n \M }_s \, < \,  + \infty.
$$
\vskip 6pt
Proof.  
\rm    In view of (3.7), it is sufficient to show that for 
     any multi-index  \f \alpha $
$$
   \sup_n \, { \M  x^{\alpha} \,  a_n \M }_0 \, < \,  + \infty.
$$
    By integration by parts, we see that
$$
x^{\alpha} \,  a_n  (x)  \, = \,   (2 \pi )^{-3} 
          \int \, e^{ix \cdot \xi}  \, 
             { \big( i \frac{\partial}{\partial \xi} \big)}^{\!\!\alpha} 
               \,  \widehat a_n (\xi) \, d\xi.               
                                                                   \T 3.8
$$
    Combining (3.8) with (3.5), (3.6) and using the
    Plancherel theorem, we get
$$
 { { \M  x^{\alpha} \, a_n \M }_0 }^2  \, \le \,
      C_{\alpha \varphi} \, 
           \int_{ \sqrt{n(n+2)} \le |\xi|  \le   \sqrt{n^2 +6n + 8} }  
             \,\,  |\xi|^{-2} \, d\xi    ,
$$
    where the constant \f C_{\alpha \varphi} \A 
    depends only on  \f \alpha \A
    and  the
    least upper bound of \f \varphi $, together with its all
    derivatives up to \f |\alpha|$--th order. This 
    gives the desired 
    \break conclusion.     $\square$ 
\vskip 12pt
     We note that the resolvent \f R_0 (z) \A  of the free Dirac 
     operator \f H_0 \A can be
     represented in 
     terms of the Fourier transform:
$$
R_0 (z) \, = \, {\SF}^{\, -1} \big[ \, 
       (  \widehat L_0 (\xi) - zI )^{-1} \, \big] \SF
            \qq ( \text{Im} \; z  \not= 0 )                    \T 3.9
$$
    where
$$
  \widehat L_0 ( \xi ) \, = \, \sum_{j=1}^3 { \xi}_j \,
         {\alpha}_j \, + \, \beta .
                                                               \T 3.10
$$
    Here an explanation must be needed. We define
     the Fourier transform of a  \f { \bC}^4$--valued function
$$
 f(x) \, = \, \pmatrix   f_1(x)   \\   f_2(x)               
       \\  f_3(x)   \\  f_4(x)   \endpmatrix  
$$
      by
$$
    \widehat f(\xi) \, = \,   
      [\SF f] (\xi) \, = \,
     \pmatrix
         \widehat f_1(\xi)    \\  \widehat f_2(\xi)
     \\    \widehat f_3(\xi)   \\   \widehat f_4(\xi)   
      \endpmatrix .
$$
     For every \f \xi \in \R3$, the Hermitian matrix
      \f \widehat L_0 (\xi)$, 
      acting in \f \bC {}^4 \A  with the usual inner product,
      has two
     eigenvalues \f \pm  \langle  \xi  \rangle $,
     each of which 
     is an eigenvalue
     of
     multiplicity two. The property
$$
  { \big( \, \widehat L_0  ( \xi )  \, \big) }^2  
      \, = \,  \langle \xi \rangle {}^2 \, I                   \T 3.11
$$
    implies that the eigenprojections 
    \f {\Psi}_{\pm} (\xi) \A 
    associated with the eigenvalues \f \pm  \langle  \xi  \rangle \A
    of \f \widehat L_0 (\xi) \A  are given by
$$
  {\Psi}_{\pm} (\xi) \, = \,  \frac12 \,
        \big( \, I \, \pm  \, \frac1{ \langle  \xi  \rangle}
           \, \widehat L_0 (\xi) \,  \big)
                                                                  \T 3.12
$$
     respectively (cf. [14, \S1]). Therefore
$$
  R_0 (z) \, = \, 
    {\SF}^{\, -1} \big[
      \, - \frac1{ \langle  \xi  \rangle  + z} \, {\Psi}_{-} (\xi)
      \, + \,
      \frac1{ \langle  \xi  \rangle  -  z} \, {\Psi}_{+} (\xi)
          \,  \big] \SF.
                                                                \T 3.13
$$
\vskip 12pt

  {\bf Lemma 3.6.} 
\SP

\it    Let \f \{ h_n \} \A be the sequence given by (3.7).   Then
$$
 \lim_{n \to \infty} \, (R_0^{\pm}(n+2)h_n, \, h_n)_0  
  \, = \, \pm  \frac{i}{4 \pi}.
$$
\vskip 6pt
Proof.  
\rm   To simplify the notation, we give the proof only for \lq\lq+".
      Since, by (3.5), \f \widehat h_n (\xi) \A  is an even
      function of \f  {\xi}_j \,  \A  and 
      \f\,  {\xi}_j  \widehat h_n (\xi) \A is an odd function
      of \f  {\xi}_j  $,
       we
      see that
$$
  \int \frac1{ \langle  \xi  \rangle  \pm z} \,
     \big\langle 
        \,   \frac1{ \langle  \xi \rangle} \,
         {\xi}_j \, {\alpha}_j \,  
            \widehat h_n (\xi),  \, \widehat h_n (\xi)   \, \big\rangle
                 \, d\xi \, = \, 0  
                  \qq  (\; j\, =\, 1, \; 2,  \; 3\;)
                                                                  \T 3.14
$$
    where  \f {\alpha}_j \A is the  matrix given in (1.1) and 
$$
    \big\langle  \,  \widehat f (\xi), \,
       \widehat g (\xi)  \, \big\rangle  \, = \,
           \sum_{k=1}^4 
             \widehat f_k (\xi) \,  
              \overline{ \widehat g_k (\xi) }.
$$
    Taking into  account (3.10), (3.12)---(3.14) and  (3.7),  we get
$$\split
  ( R_0 (z) h_n, \, h_n)_0 \, & = \,  (2 \pi )^{-3} \, 
       \int \frac{-1}{ \langle  \xi  \rangle  +  z} \,
         \cdot
         \frac12 (1 \, - \, \frac1{ \langle \xi \rangle}) \,
         {\big|  \,  \widehat a_n (\xi) \big|}^2         \, d\xi    \\
    & \qq + \;
          (2 \pi )^{-3} \, 
       \int \frac{1}{ \langle  \xi  \rangle  -  z} \,
         \cdot
         \frac12 (1 \, + \, \frac1{ \langle \xi \rangle}) \,
         {\big|  \,  \widehat a_n (\xi) \big|}^2         \, d\xi , 
\endsplit
$$
     where we used (3.1).  Using (3.5) and 
     passing to the polar coordinates, we have
$$
\left\{
\aligned
     {} &  4   \pi {}^2 \,  ( R_0 (z) h_n, \, h_n)_0 \, \\ 
       & = \, - \, \int_0^{\infty}  \,
          \frac1{\sqrt{r^2 + 1} + z} \, 
           \big( 1 -  \frac1{ \sqrt{r^2 + 1} }  \big)  \,
              \varphi( \sqrt{r^2 + 1}  - n -2) ^2 \, dr     \\
       & \ \   +    \int_0^{\infty}  \,
          \frac1{\sqrt{r^2 + 1} -  z} \, 
           \big( 1 +  \frac1{ \sqrt{r^2 + 1} }  \big)  \,
              \varphi( \sqrt{r^2 + 1}  - n -2) ^2 \, dr        \\
       & =   \, - \,  \int_{-1}^1
          \,  \frac1{ \sigma + n + 2 + z}  
            \,  \big( 1 -  \frac1{ \sigma + n +2 }  \big) \,  
               \varphi ( \sigma )^2  \,
                \frac{\sigma + n + 2}{  \sqrt{(\sigma + n + 2)^2  -1}   }
                      \,   d\sigma                                     \\
      & \ \   +     \int_{-1}^1
          \,  \frac1{ \sigma + n + 2 - z}  
            \,  \big( 1 +  \frac1{ \sigma + n +2 }  \big) \,  
               \varphi ( \sigma )^2  \,
                \frac{\sigma + n + 2}{  \sqrt{(\sigma + n + 2)^2  -1}   }
                      \,   d\sigma.              
\endaligned  \right.                                               \T 3.15
$$
     In the second equality above, we made a change of a variable.
     Putting 
$$
    z \, = \, n + 2 + i \mu  \qq  ( \mu > 0),
$$
     and taking the limit of (3.15) as  $\,  \mu  \downarrow 0$,
     we see, 
     in view of Theorem 2.2 and Lemma 3.3, that
$$
\left\{
\aligned
      4   \pi {}^2 \, &  ( R_0^+ (n + 2) h_n, \, h_n)_0 \, \\ 
       = &  \, - \,  \int_{-1}^1
          \,  \frac1{ \sigma + 2 n + 4 }  
            \,  \big( 1 -  \frac1{ \sigma + n +2 }  \big) \,  
               \varphi ( \sigma )^2  \,
                \frac{\sigma + n + 2}{  \sqrt{(\sigma + n + 2)^2  -1}   }
                      \,   d\sigma                                     \\
      {}  &   \q  +  \,   i \pi \, ( 1 + \frac1{n + 2} ) \,
           {\varphi (0)}^2  \,  
                 \frac{n + 2}{ \sqrt{(n+2)^2 -1} }                   \\
      {}   &   \qq +  \,  
       \int_{-1}^1 \big\{ \int_0^1   { \omega }_n^{\prime}
          ( \sigma  \theta )  \,  d\theta  \,    \big\} \, d\sigma, 
\endaligned    \right.                                       
                                                                  \T 3.16
$$
     where 
$$
   { \omega }_n (\sigma) \, =\, 
    \big( 1 +  \frac1{ \sigma + n +2 }  \big) \,  
               \varphi ( \sigma )^2  \,
                \frac{\sigma + n + 2}{  \sqrt{(\sigma + n + 2)^2  -1}   }.
$$
   It is easy to see that the integrand in  the first term on the 
   right hand side of (3.16)  is
   less than or equal to,
   in the absolute value, \f K/n  \A  where \f K \A   is a 
   positive constant independent of
   \f n$.  Therefore, the first term  converges to  0 as 
   \f n \to  \infty$. 
   As for the third term on the 
   right hand side of (3.16), we see that  \f \{ { \omega}_n^{\prime}
    \} _{n=1}^{\infty}  \A  is 
    an 
    uniformly bounded sequence of functions
    which converges  pointwisely  
     to  \f  ( { \varphi }^2  ) ^{\prime} $.
    Hence, by the
    Lebesgue dominated convergence theorem, 
$$
  \lim_{n \to \infty} \, 
 \int_{-1}^1 \big\{ \int_0^1   { \omega }_n^{\prime}
          ( \sigma  \theta )  \,  d\theta  \,    \big\} \, d\sigma
  \,   =   \,
   \int_{-1}^1 \big\{ \int_0^1   {( \varphi {}^2) }^{\prime}
          ( \sigma  \theta )  \,  d\theta  \,    \big\} \, d\sigma.
                                                                \T 3.17
$$
    In view of the fact that  \f {\varphi}^2 \A is an 
    even function, it follows
    from Lemma 3.4 that the  right hand side of (3.17)
    equals 0. 
    Summing up, we get
$$
  \lim_{n \to \infty}
  4   \pi {}^2 \,  ( R_0^+ (n + 2) h_n, \, h_n)_0 \,
   = \, i \,  \pi  \,  \varphi (0) {}^2  
     \, = \,  i \,  \pi  .
$$
    This completes the proof.  \f  \square$
\vskip 12pt
    It is obvious that Lemmas 3.5 and 3.6  give the proof of 
    Proposition 3.1.


\vskip 24pt  

\heading      {\bf \S4. A known result for 
               Schr\"odinger operators }              \endheading 

\SP

     The limiting absorption  principle for Schr\"odinger operators
     has been extensively studied in connection
     with the spectral
     and scattering theory; cf. [8], [3], [1].
     In this section, we make a brief review of
     a result due to Sait\B o\;[9] and [10], which will
     be used 
     in the proof of Theorem 2.4. 
     Let \f T \A denote the selfadjoint 
     operator which is defined to be the closure of \f  - \Delta +  V(x) \A
     restricted on \f \Con(\Rn)$, where \f V(x) \A is  a real-valued
     function satisfying
$$
          | \,  V(x)\, |  \, \le \, C \,
                              { \langle x \rangle }^{-1-\epsilon}
                                                                  \T 4.1
$$
     for  $\, C > 0 \, $ and $\, \epsilon > 0  \,$. Let 
$$
    \Gamma (z) \, 
     = \, (T \, - \, z)^{-1}.
$$
     Then it is well-known that the limiting 
     absorption principle holds for $\, T \,$, 
     \break that is, for any   
     $\,  \lambda > 0 \,$, there correspond the extended resolvents
     $\,  { \Gamma }^{\pm}(\lambda) \,$ in
     \break $\bB(L_{2,s}(\Rn), L_{2, -s}(\Rn)) $ such that  
     for any \f f \A in   \f L_{2, s}( \Rn )$
$$
   \Gamma  ( \lambda \pm i \mu ) f  \,  \longrightarrow  \, 
        {\Gamma}^{\pm} (\lambda) f 
     \q \hbox{ in } \ L_{2, -s}
                                                       \T 4.2
$$
    as \f  \mu \downarrow 0 $. Furthermore, it is known that  
     \f {\Gamma}^{\pm}(\lambda)f \A  are $L_{2, -s}(\Rn)$--valued 
     continuous functions in \f (0, \, \infty)$. (Sait$\bar {\hbox{o} }$
     \;[8], Ikebe-Sait$\bar {\hbox{o} }$\;[3] and Agmon\;[1].)
     As for asymptotic behaviors of ${\Gamma}^{\pm} (\lambda)$, we have

\vskip 12pt
        { \bf Theorem 4.1 } {\it (Sait$\bar {\hbox{o} }$\;[9, 10]).
$$
   {  \M  {\Gamma}^{\pm} ( \lambda )  \M }_{(s, -s)}
          \, = \,  O( {\lambda}^{-1/2})
                              \qq ( \lambda  \to  \infty ).    \T 4.3
$$
\vskip 12pt

\rm      Also, Sait$\bar {\hbox{o} }$  proved

\vskip 12pt

        {\bf Theorem 4.2}\ {\it (Sait$\bar {\hbox{o} }$\;[9, 10]).
\SP 

        Let  \f s >  1/2$.  Then for any  \f d \, > 0 \A there exists 
     a positive constant \f C > 0 \a such that 
$$
     {  \M {\Gamma} ({ \kappa }^2 ) \M }_{(s, -s)}
                        \,  \le  \, C / |\kappa|
                                                          \T 4.4
$$
     for all \f \kappa \A  with \f | \hbox{\rm Re} \, \kappa |  > d \, \A
     and \f \hbox{\rm Im} \, \kappa > 0 $.}
\vskip 12pt
       {\bf Remark 4.3}
\SP      
       
\rm      We would like to emphasize that the asymptotic
        behavior (4.3) is 
       useful in the inverse scattering theory for the Schr\"odinger
       operator. See
        Sait$\bar {\hbox{o} }$\;[11, 12].


\vskip 24pt

\heading      {\bf \S5. Pseudodifferential operators }   \endheading 

\SP

        The proof of Theorem 2.4 is based on the resolvent 
        estimate for the
      Schr\"odinger operator (Theorem 4.2) 
      as well as the theory of 
        pseudodifferential operators. In this section we  introduce 
     a class of
     symbols of pseudo-differential operators which are suitable to our 
     purpose.  We then  establish  boundedness results
     in the weighted Hilbert spaces, which are important in
     relation to
     the limiting absorption principle for the Dirac operator;
     cf. [7].

\vskip 12pt

        {\bf Definition 5.1.}
\SP 

        A \f C^{\infty}\A  function \f p(x, \xi) \A
         on \f \R3\times\R3 \A
     is said to be in the class \f S_{0,0}^{m} \;\; ( m \in \Ro) \A 
      if for any pair 
     \f \alpha \A  and \f \beta \A 
      of multi-indicies there exists a constant
     \f C_{\alpha\beta} \ge 0 \A  such that
$$
       \big| \big(\frac{\pa}{\pa\xi}\big)^{\alpha}                    
           \big(\frac{\pa}{\pa x}\big)^{\beta} p(x, \xi) \big|
                                        \le C_{\alpha\beta}      
                             \, \langle \xi \rangle {}^{\! m} 
$$
     for all  \f x, \xi \in \R3 $.

\vskip 12pt

        {\bf Remark 5.2.}
\SP 

        The class \f S_{0,0}^{m} \A  is a Fr\'echet space equipped with 
     the semi-norms
$$
    |p|_{\ell}^{(m)} = \max_{|\alpha|,|\beta|\le\ell} \sup_{x,\xi}\bigg\{                                                   
          \,  \langle  \xi  \rangle {}^{\! -m}
          \big| \big(\frac{\pa}{\pa\xi}\big)^{\alpha}                    
           \big(\frac{\pa}{\pa x}\big)^{\beta} p(x, \xi) \big|\bigg\}
                       \q (\ell = 0, 1, 2, \cdots).        
$$
\vskip 12pt

        A pseudodifferential operator $\,p(x, D)\;$ with symbol 
     \f p(x, \xi) \A  is defined by
$$
      p(x, D)f(x) = (2\pi)^{-3} \int_{\R3} e^{ix\cdot\xi}\, p(x, \xi)
                      \,   \widehat f(\xi)\, d\xi                  
$$
     for \f f \in \CS(\R3)$.

\vskip 12pt 

        {\bf Lemma 5.3.}
\SP 
 
        \it  Let \f p(x, \, \xi) \A  be in  \f S_{0,0}^0 \A 
        and let \f s > 0$.   Define the oscillatory integral
$$
   r(x, \, \xi) \, = \,
     \text{\rm Os}-\! \iint_{{\bold R}^6} \;
        e^{-iy \cdot  \eta}   \; p(x, \,  \xi + \eta) \;
            \langle x + y \rangle {}^{\! - s} \;  
               (2 \pi )^{-3}\,  dyd\eta  .
$$
        Then for any pair \f \alpha \A and \f \beta \A  
        there exists a
        constant \f C_{s \alpha \beta}  \ge  0 \A
        such that 
$$
       \big| \big(\frac{\pa}{\pa\xi}\big)^{\alpha}                    
           \big(\frac{\pa}{\pa x}\big)^{\beta} r(x, \xi) \big|
                                        \le C_{s\alpha\beta} 
        \,    |p|_k^{(0)}  \,
             \langle x  \rangle {}^{\! -s},  
$$
        where
$$
   k \, = \, \max \,  \big\{ \,
        | \beta |, \; 2[s + \frac52] + |\alpha| \,  \big\}.
                                                             \T 5.1
$$
\rm
\vskip 12pt
    {\bf Remark 5.4.}
\SP

   \item{(i)}  For the definition of oscillatory integral,
     see Kumano-go[5, Chapter 1,
     Section 6\,].
   \item{(ii)} For a positive number \f s$, [$s$]  denotes 
   the largest integer less than or equal to \f s $.
   \item{(iii)} By [5, Theorem 2.6(1), p. 74], \f r(x, \, D) \, =
    \, p(x, \, D) \, \langle x \rangle {}^{\! -s}$.
    This fact will be used in the proof of Lemma 5.5 below.

\vskip 15pt
\it Proof.
\rm     By differentiation under the oscillatory integral sign\,
        (cf. [5, (2.23),
        p. 70]), we see that
$$\aligned
  {}& \big(\frac{\pa}{\pa\xi}\big)^{\alpha}                    
           \big(\frac{\pa}{\pa x}\big)^{\beta} r(x, \xi)      \\
   {} &=\,  
       \text{Os}-\! \iint_{{\bold R}^6} \;
        e^{-iy \cdot  \eta}   \;
        \big(\frac{\pa}{\pa\xi}\big)^{\alpha}                    
           \big(\frac{\pa}{\pa x}\big)^{\beta}
           \big\{ p(x, \,  \xi + \eta) \;
            \langle x + y \rangle {}^{\! - s}  \big\}
              \;   (2 \pi )^{-3}   dyd\eta.  
       \endaligned                                      \T 5.2
$$
       Putting
$$
     M  \, = \, \big[ \frac{ s + 5}2 \big] ,
$$
     and   integrating by parts,  we see that 
$$\aligned
   \text{  RHS of (5.2)}  \, = \, 
       \iint_{{\bold R }^6 }   \; & e^{-iy \cdot  \eta} \;
         \langle y \rangle {}^{\! -2M} \,
               \langle D_{\eta} \rangle {}^{\! 2M} \,
         \big\{  
             \langle \eta  \rangle {}^{\! -4} \,
             \langle D_y \rangle {}^{\! 4} \,     \times   \\
       \times & 
        \big(\frac{\pa}{\pa\xi}\big)^{\alpha}                    
           \big(\frac{\pa}{\pa x}\big)^{\beta}
           \big\{ p(x, \,  \xi + \eta) \;
            \langle x + y \rangle {}^{\! - s}  \big\}
              \;   (2 \pi )^{-3}   dyd\eta ,
\endaligned                                               \T 5.3    
$$
       where 
$$
      \langle D_{\eta}  \rangle {}^2  \,
           = \, 1 \, - \, { \Delta}_{\eta} \;  ;
      \q
      \langle D_y  \rangle  {}^2 \,  = \, 1 \, - \, { \Delta}_y 
$$
     (cf. [5, Theorem 6.4, p. 47]).  Note that the integral
     on the right
     hand side of (5.3) is in the usual sense.
     Using the following two inequalities
$$
   \big|   \big(\frac{\pa}{\pa x}\big)^{\alpha} 
        \,   \langle  x   \rangle {}^{ \! -s}   \big|  \,
        \le \,
          C_{\alpha s} \, \langle x  \rangle {}^{\! -s}
$$
      and
$$
      \langle  x + y  \rangle {}^{ \! -1}  \,   \le  \,
         \sqrt2     \langle y  \rangle   \langle x  \rangle {}^{\! -1},
$$
     we get 
$$\aligned
    {} &  |\;  \text{the integrand on the RHS of (5.3) }  \, |    \\
     {} & \qq  \le  \,  C_{s\alpha \beta} \,
         |p|_k^{(0)}  \,      \langle x  \rangle {}^{\! -s}
          \,  \langle y  \rangle {}^{\! -2M+s}
          \,      \langle \eta  \rangle {}^{\! -4} ,
\endaligned                                                    \T 5.4  
$$
      where \f C_{s \alpha \beta} \A  is a nonnegative constant.
      We note that
      \f     \langle y  \rangle {}^{\! -2M+s}
          \,      \langle \eta  \rangle {}^{\! -4} \A
      is integrable on \f { \bold R}^6$.  
      Hence, combining (5.2)--(5.4), we get the desired conclusion.
                                                            \f \square$  
\vskip 12pt

       {\bf Lemma 5.5}
\SP

       {\it Let \f p(x, \xi) \A  be in \f S_{0,0}^{0} $. Then for any 
     \f s \ge  0 \A  there exist a nonnegative  constant \f C \A and
     a positive integer \f \ell  \A such that }
$$
       \M p(x, D)f \M_s \le C \, |p|_{\ell}^{(0)} \, \M f \M_s
                            \qq (f \in \CS(\R3) \,) ,             \T 5.5
$$
    {\it  where \f C \A and \f \ell \A  depend only on \f s$.  }
\vskip 6pt
\it    Proof. \rm
      It is sufficient to show that for any \f s > 0 \A  there exist
      a nonnegative  constant \f C \A and  a positive 
      integer \f \ell \,  \A such that
$$
      \M \,    \langle x  \rangle {}^{\! s}  \, 
        p(x, \, D)  \,   \langle x  \rangle {}^{\! - s} \, f 
           \,  \M {}_0
         \;  \le  \;  
          C \,  |p|^{(0)}_{\ell} \,
           \M f \M {}_0                                       \T 5.6
$$
      for all \f f \in \Cal S (\R3)$.
      Let \f r(x, \, \xi) \A be the symbol defined in Lemma 5.3,
      and put
$$
     q( x, \, \xi) \, = \,  \langle x  \rangle {}^{\! s}  \, 
                    r(x, \, \xi).
$$
    According to Remark 5.4(iii),
$$
     q( x, \, D) \, = \,  \langle x  \rangle {}^{\! s}  \, 
                    p(x, \, D) \,  \langle x  \rangle {}^{\! -s} .
$$
     It follows from Lemma 5.3 that for any pair \f \alpha \A
     and \f \beta \A  there exists a constant
     \f C_{s\alpha\beta}  \ge  0 \A  such 
     that
$$
  \big|
     \big(\frac{\pa}{\pa\xi}\big)^{\alpha}                    
           \big(\frac{\pa}{\pa x}\big)^{\beta} q (x, \xi)   \big|
    \;  \le \;    C_{s\alpha\beta}  \,
          |p|^{(0)}_k
$$
        where  \f k \A  is given by (5.1).
        Then the Calder\'on-Vaillancourt theorem([2], [5,
        Theorem 1.6, p.224]) implies (5.6).
                                                \f \square$.
%
%
%

\newpage
        {\bf Lemma 5.6.}
\SP

  {\it Let \f p( \xi) \A  be in  \f S_{0,0}^{-1} $. Then for any 
     \f s  \ge  0 \A  there exist a nonnegative  constant \f C \A and
     a positive integer \f \ell  \A such that }
$$
       \M p( D)f \M_{1,s} \le C \, |p|_{\ell}^{(-1)} \, \M f \M_s
                            \qq (f \in \CS(\R3) \,) ,             \T 5.7
$$
    {\it  where \f C \A and \f \ell \A  depend only on \f s$.  }
\vskip 6pt
\it    Proof. \rm
    By definition (1.9), we have
$$
   \M p( D)f {\M_{1,s}}^{\!\!\!\!  2} 
      \, = \,   \M p( D)f { \M_s}^{\!\! 2} 
      \, + \,  \sum_{j=1}^3  \,
           \M  \frac{\partial}{\partial x_j} p( D)f { \M_s }^{\!\! 2}.
                                                                 \T 5.8 
$$
    Regarding \f p(\xi) \A as a symbol in \f S_{0,0}^0 $,
    and applying Lemma 5.5, we get
$$
  \M p( D)f \M_s   \, \le \,  C \, |p|_{\ell}^{(0)} \, \M f \M_s
                                 \qq (f \in \CS(\R3) \,).          \T 5.9
$$
    Note that the symbol of 
    \f  (\partial / \partial x_j ) p( D) \A  is \f i \,\xi_j \,  p(\xi)$,
    which belongs to \f S_{0,0}^0 $.  Then, by Lemma 5.5, 
    we see that
$$
\M  \frac{\partial}{\partial x_j} p( D)f \M_s  \, \le \, 
       C \,  |p|_{\ell}^{(-1)} \, \M f \M_s
                                 \qq (f \in \CS(\R3) \,).        \T 5.10
$$
    Using the fact that \f  |p|_{\ell}^{(0)}  \,  \le \,
      |p|_{\ell}^{(-1)} \A for \f \ell = 0, \, 1, \, 2, \cdots $,
   and combining  (5.8)--(5.10), we obtain (5.7).
                                                              $\square$
\vskip 12pt
     We now  need to extend Lemmas 5.5 and 5.6
     to a system of pseudodifferential operators. Let 
$$
         P(x, \xi) ={ \big( p_{jk}(x, \xi) \big)}_{1\le j,k \le 4}
$$
     be a \f 4 \times 4 \A  matrix-valued symbol. Then we define
$$
         P(x, D ) = {\big( p_{jk}(x, D ) \big)}_{1\le j,k \le 4}
$$
     by 
$$
      P(x, D)f(x) = (2\pi)^{-3} \int_{\R3} e^{ix\cdot\xi} \, P(x, \xi)
                       \,   \widehat f(\xi)\, d\xi                 
$$
     for \f f \in [ \CS(\R3)]^4$.  
          If  \f  p_{jk}(x, \xi )
            \in  S_{0,0}^{\,m\,}$, \f 1\le j,k \le 4 $,
     we define
$$
      |P|_{\ell}^{(m)} =   \big \{  \sum_{j, k=1}^4 
                   ( |p_{jk}|_{\ell}^{(m)} )^2 \,   \big\} ^{1/2}   \T 5.11        
$$
     for \f \ell = 0, 1, 2, \dots $, where \f  |p_{jk}|_{\ell}^{(m)}  \A  
     are the semi-norms introduced in Remark 5.2. 
     We then have  natural extensions of Lemmas  5.5 and  5.6.
%
%
%

\newpage
        {\bf Lemma 5.7.}
\SP

        {\it  Let  \f  p_{jk}(x, \xi )  \A be in
                    \f  S_{0,0}^{\,0\,} \A
     for  \f j, \, k = 1, \, 2, \, 3, \, 4$. 
     Then for any \f s \ge  0 \A  there 
     exist a nonnegative  constant \f C \A  and a positive integer
     \f \ell \A such that 
$$
     { \M  P(x, D ) f \M }_s \,  \le \,
               C \,  |P|_{\ell}^{(0)} \, { \M f \M }_s  
                \qq (\, f \in  [\CS(\R3)]^4 \,) ,  
                                                            \T 5.12    
$$
      {\it   where \f C \A and \f \ell \A  depend only on \f s$.}
\vskip 6pt
\it      Proof.   \rm
         It is a matter of simple  computation:
$$
\aligned
   {  \M  P(x, \, D) f \M {}_s  }^{\! 2}  \, = \, 
        \sum_{j=1}^4  \;     &
         \M   \sum_{k=1}^4  \, 
              p_{jk}(x, \, D) \, f_k  \, { \M_s}^{\! 2}      \\
   \le \,    \sum_{j=1}^4  \, 
       \big(    \,   \sum_{k=1}^4  \, & C \, |p_{jk}|^{(0)}_{\ell}  \,
          \M f_k \M_s  { \big)}^2  
              \q   ( \text{by Lemma 5.5} )                       \\
    \le \, C^2    
          \, \sum_{j=1}^4 \, \big\{ \,  \sum_{k=1}^4  \,
                 &  { \big( \, |p_{jk}|_{\ell}^{(0)}  \, 
                      \big) }^2  \,
                \big\} \,  \sum_{k=1}^4 
                  \,  { \M f_k \M_s}^{\! 2} 
                    \q  ( \text{by  the Schwarz inequality} ).
\endaligned    
$$
   With the notation (5.11), this is equivalent to (5.12).       \f \square$
\vskip 12pt
        {\bf Lemma 5.8.}
\SP

        {\it  Let  \f  p_{jk}(\xi )  \A be in
                    \f  S_{0,0}^{\,-1\,} \A
     for  \f j, \, k = 1, \, 2, \, 3, \, 4$. 
     Then for any \f s \ge  0 \A  there 
     exist a nonnegative  constant \f C \A  and a positive integer
     \f \ell \A such that 
$$
     { \M  P( D ) f \M }_{1,s} \,  \le \,
               C \,  |P|_{\ell}^{(-1)} \, { \M f \M }_s  
                \qq (\, f \in  [\CS(\R3)]^4 \,) ,  
$$
      {\it   where \f C \A and \f \ell \A  depend only on \f s$.}
\vskip 12pt
\rm
    In view of Lemma 5.6, the proof of Lemma 5.8 is 
    similar to that of Lemma 5.7. We should like to 
    mention that Lemma 5.8 is beyond the necessity
    for the present paper. However, we need the lemma
    in our forthcoming paper\,[7].


\vskip 24pt

\heading      {\bf \S6. Proof of Theorem 2.4  }   \endheading 

\SP

    In this section, we give the proof of Theorem 2.4. 
    We begin with rewriting (3.9). 
    Using (3.11), we see that
$$
 ( \widehat L_0 (\xi ) \, - \, zI ) \,
         ( \widehat L_0 (\xi ) \, + \, zI )  \; = \; 
            ( \, \langle \xi \rangle {}^2 \, - \, z^2 \, ) \, I.
                                                                 \T 6.1
$$
      Hence
$$
R_0 (z) \, = \, {\SF}^{\,-1} \big[ \, 
     \frac1{ \langle \xi \rangle {}^2 \, - \, z^2 } 
         \, (  \widehat L_0 (\xi) + zI )
             \, \big] \SF
            \qq ( \text{Im} \; z  \not= 0 ).                    \T 6.2
$$

\vskip 12pt

     {\bf Theorem 6.1.}
\SP

     \it   Suppose that \f s > 1/2$. Then
$$
  \sup \, \big\{ \; 
       \M R_0 (z) \M_{(s, -s)} \;  \big/  \;\;\;
          2 \, \le \, | \text{\rm Re} \; z |, \;\;  0 \, < \, 
             | \text{\rm Im} \; z | \, < \, 1  \;  \big\}  
               \; < \; +\infty.                              \T 6.3
$$
 \rm

\vskip 12pt

      {\bf Remark 6.2.}
\SP

      It is evident that Theorem 6.1, together with Theorem 2.2, 
      implies Theorem 2.4.

\vskip 12pt

       {\it  Proof.}  \  Set 
$$
     J = \big\{  \,  z \in  \bC \, \big/ 
      \;\;\;
          2 \, \le \, | \text{Re} \; z |, \;\;  0 \, < \, 
             | \text{Im} \; z | \, < \, 1  \;  \big\}.
$$
     Choose  \f \rho  \in  C_0^{\infty}( {\bold R}) \A  so that
$$
      \rho (t) = \cases  1, & \text{ if \ $|t| <1/2$ }   \\  
            {} &  \text{}                              \\
                  0, &  \text{ if \ $ |t| > 1$ }  .  
     \endcases              
$$
     For each \f z \in J$, we define a cutoff function 
     \f { \gamma}_z(\xi) \A on  \f \R3 \A by 
$$
      { \gamma}_z(\xi) = \cases  \rho(
       \langle  \xi  \rangle  - \hbox{Re}   \; z), &
         \text{  if \ Re$\;z \ge 2$ }                          \\
         {} & \text{}                                     \\
       \rho(
        \langle  \xi  \rangle   + \hbox{Re}   \; z), &
          \text{  if \ Re$\;z \le -2$} .
      \endcases   
$$
     Using (6.2) and \f { \gamma}_z(\xi) $, we decompose the resolvent of
     \f H_0 \A into three parts:
$$
        R_0(z) = ( -\Delta  + 1 -z^2)^{-1}\, A_z \, + \, B_z 
                        \, + z( -\Delta  + 1 -z^2)^{-1}
$$
     where
$$
\align
       A_z &=  {\SF}^{\,-1}  \Big[ {\gamma}_z (\xi) { \widehat L}_0  (\xi)
          \Big]  {\SF},                                   
\\
       B_z &= {\SF}^{\,-1}  \Big[ 
         \frac{ 1 - {\gamma}_z (\xi) }
            { { \langle  \xi  \rangle  }^2 - z^2  }
                         { \widehat L}_0  (\xi)  \Big]  {\SF}.
\endalign
$$
    Note that for \f  \xi \, \in \,  \text{supp}[ \, \gamma_z\,] \A
    with \f z \in J$
$$
    \frac14 \, |z| \, \le \, 
         \langle  \xi  \rangle   \,   \le  
           \frac32 \, |z|.                              \T 6.4
$$
    Using (6.4) and (3.10), we see that for any  \f \alpha \,\A there
    exists a constant \f C_{\alpha} \A 
    such that
$$
    \big|  \,  { \big( 
        \frac{\partial}{\partial \xi}  \big) }^{\! \alpha}
         \big(  \,  \gamma_z (\xi) \,
            \widehat L_0 (\xi)  \, \big)  \, \big|
            \,  \le  \,
             C_{\alpha}  \, |z|                              \T 6.5
$$
    for all \f z \in J$.
    Here and in the sequel, for a \f 4 \times 4 \A matrix \f M$, 
    its matrix norm is denoted by \f |M| \A (e.g.,
    \f |M|^2 \, = \, \sum_{j, k=1}^4 \, {m_{jk}}^2\A; 
    actually it is irrelevant  
    which norm one chooses), and
    for a \f 4 \times 4 \A matrix-valued function
    \f M(\xi) \, = \, \big( \, m_{jk}(\xi) \, \big)$
    we write
$$
  { \big( 
        \frac{\partial}{\partial \xi}  \big) }^{\! \alpha} \, M(\xi) \,
        =
  \Big(  
      { \big( 
        \frac{\partial}{\partial \xi}  \big) }^{\! \alpha} \,
            m_{jk}(\xi)   \,    \Big)_{1 \le j, \, k  \le 4}.
$$
   Then noting (6.5) and
     applying Lemma 5.7  to \f A_z$, we get
$$
      { \M  A_z f \M }_s \,  \le \,  C_1 |z| \,  {\M f \M}_s
            \qq (\, f \in  [\CS(\R3)]^4 \,)
                                                      \T 6.6
$$
    for \f z \in J$,  
     where \f C_1 \A  is independent of \f z \in J$. 
     On the other hand, it follows,
     in particular, from Theorem 4.2 that for \f z \in J$
$$
{ \M ( -\Delta  + 1 -z^2)^{-1} \M  }_{(s, -s)}  \, \le  \, 
   \frac{ \, C_2 \, }{  |z|}                                   \T 6.7
$$
    with a constant \f C_2 \A independent of \f z$. 
     Combining (6.6) and (6.7), we have
$$
        { \M ( -\Delta  + 1 -z^2)^{-1} A_z f \M  }_{-s}  \,
               \le \,  C_3  {\M f \M}_s 
                \qq (\, f \in  [\CS(\R3)]^4 \,) ,
                                                     \T 6.8
$$
     where  \f C_3 \A  is independent of \f z \in J$.

     In order to apply Lemma 5.7 to \f B_z$, we note that
$$
  \big| \, \langle \xi \rangle {}^2  \, - \, z^2 \, \big|
     \, \ge \,  \frac12 \, \langle \xi \rangle
                                                          \T 6.9
$$
   for \f z \in J \A and \f \xi \in  \text{supp}[ 1- \gamma_z ]$.
   Using (6.9), we see that
   for any \f \alpha \A there exists a 
   constant  \f C_{\alpha}^{\prime} \A  such that
$$
  \big|   { \big( 
        \frac{\partial}{\partial \xi}  \big) }^{\! \alpha} 
        \big\{ \frac{ 1 - {\gamma}_z (\xi) }
            { { \langle  \xi  \rangle  }^2 - z^2  }  \big\}   \big|
      \;  \le  \; 
         C_{\alpha}^{\prime} \, 
           \langle \xi \rangle {}^{ -1}
$$
    for all \f z \in J \A and all \f \xi \in \R3$.
    Hence, to each \f \alpha \,\A we can find
     a constant \f C_{\alpha}^{\prime\prime} \A satisfying
$$      
  \big|   \,  { \big( 
        \frac{\partial}{\partial \xi}  \big) }^{\! \alpha} 
        \big[   \frac{ 1 - {\gamma}_z (\xi) }
            { { \langle  \xi  \rangle  }^2 - z^2  }
              \,  \widehat L_0 ( \xi ) \,  \big]  \, \big|
                \, \le  \, C_{\alpha}^{\prime\prime}
$$
    for all \f z \in J\A and all  \f  \xi \in  \R3$. 
    We then apply Lemma 5.7 to \f B_z$, and 
    deduce that
$$
      { \M  B_z f \M }_{-s} \,  \le \,  C_4 {\M f \M}_s
           \qq (\, f \in  [\CS(\R3)]^4 \,)
                                                               \T 6.10
$$
   with a constant \f C_4\A independent of \f z \in J$. 
   Since   \f [\CS(\R3)]^4  \A   is dense in   
    \f \CL_{2,s}$, we conclude from (6.8), (6.10) and (6.7) that
    (6.3) holds.   $\square$


\vskip 24pt

\heading      {\bf \S7. Proof of Theorem 2.5 }   \endheading 

\SP

    In order to prove Theorem 2.5, we need some prerequisites and
    a few lemmas. Throughout this section, we
    regard  \f \CS(\R3) \A  as a Fr\'echet space
    equipped with the semi-norms
$$
       |a|_{\ell , \, \CS} \,  = \,
      \sum_{ | \alpha + \beta | \le  \ell}  \;
          \sup_x \, \big\{ \, \big| \,
             x^{\alpha} \, 
     { \big(  \frac{\partial}{\partial x}  \big) }^{\! \beta} 
                                          a(x) \big|   \, \big\}   
        \qq   ( \ell \, = \, 0, \, 1, \, 2, \, \cdots \, ).
                                                                   \T 7.1
$$
    For \f f \in [\CS(\R3)]^4  \, \A  we introduce semi-norms
    by
$$
 |f|_{\ell , \, \CS} \,  = \,
      \sum_{k=1}^4  \, |f_k|_{\ell , \, \CS} 
            \qq ( \ell \, = \, 0, \, 1, \, 2, \, \cdots \, ).
$$
    It is  then trivial that   \f  [\CS(\R3)]^4 \A   is a 
    Fr\'echet space. Note that
    we use the same notation  \f  | \cdot |_{\ell , \, \CS} \A  
    as in (7.1).  We believe that 
    this causes no confusion.  For                           \break
    \f f \in [\CS(\R3)]^4  \,\A  we define
$$
   \text{supp}[\,f\,] \, = \, \bigcup_{k=1}^4 \,
                 \text{supp}[\,f_k\,] .
$$
\vskip 12pt
 
      {\bf Lemma 7.1.}
\SP

        {\it Define}
$$
     { \CX}_0  \, = \, \big\{  \,  f \in [ \CS({\bold R}_x^3)]^4 
         \, \big/  \,\,   
            \CF f  \in  \,[\Con({\bold R}_{\xi}^3)]^4\,   \big\}.          
$$
     { \it  Then \f { \CX}_0  \a   is dense in  \f  \CL_{2,s} \A
     for any  \f s \in {\bold R} $ }. 
\vskip 6pt
    {\it Proof.}   Let \f s \A be in \f {\bold R}$.  Let 
    \f g  \in    \CL_{2,s} \A  and
    \f  \epsilon > 0 \A be given.
    Since \f  [ \CS(\R3)]^4 \A  is dense in    \f  \CL_{2,s} $,
    we can find  \f f_{\epsilon} \in [ \CS(\R3)]^4 \A
    such that
$$
    \M \, g \, - \, f_{\epsilon} \, \M_s   \;  <  \;
            \frac{\epsilon}2.                               \T 7.2
$$
    Note that \f [ C_0^{\infty}( \R3 ) ]^4 \A  is
     dense in    \f [ \CS(\R3)]^4 $.  We then see that there
     exists a sequence \f \{ v_n \}_{n=1}^{\infty} \,
     \subset \, [\Con({\bold R}_{\xi}^3)]^4  \A
     such that
$$
  v_n  \, \longrightarrow \, \SF f_{\epsilon}  \;\;
    \text{ in }\;   [ \CS({\bold R}_{\xi}^3 )]^4  \;
       \text{ as } \; n \to \infty.
                                                              \T 7.3  
$$
   Now put
$$
   g_n \, = \, {\SF}^{\, -1} v_n  \q  ( n = 1, \, 2, \,  \cdots ).
                                                                \T 7.4
$$
    Since \f {\SF}^{\, -1}  \A  is a continuous map from
    \f  [ \CS({\bold R}_{\xi}^3 )]^4  \A to
     \f  [ \CS({\bold R}_x^3 )]^4 $, we deduce from
    (7.3) that \f g_n \to  f_{\epsilon} \A in
     \f  [ \CS({\bold R}_x^3 )]^4  \A
     as \f n \to \infty$. In particular, we have
$$
  \M \, g_n \, - \, f_{\epsilon} \, \M_s   \to 0 
                       \q ( n \to \infty).
$$
   Therefore we can choose an integer \f N \A so that
$$
\M \,  f_{\epsilon} \, - \,  g_N \,  \M_s  \; < \;
        \frac{\epsilon}2.                                     \T 7.5  
$$
    Then we see, by (7.2) and (7.5), that
$$
\M \,  g  \, - \,  g_N \,  \M_s  \; < \;   \epsilon,
$$
    and, by (7.4), that \f g_N \in  {\CX}_0$.             $\square$
\vskip 12pt
     
        {\bf Lemma 7.2.}
\SP

        {\it  For \f z \in \bC $, put }
$$
       R(\xi; z) \, = \, 
            \frac1{ { \langle  \xi  \rangle  }^2 - z^2 } \,
       (  {\widehat L}_0  (\xi) +  zI ) .
$$
   {\it  Then for any \f K > 1 \A  and any multi-index  \f \alpha \A
     there exists a constant  \f C_{\alpha K} >0 \A such that}
$$
    \big| \big( \frac{\pa}{\pa\xi} \big)^{\! \alpha} 
            \, R(\xi; z) \big|
      \,  \le  \, \frac{ \, C_{\alpha K} \, }{ |z| }         
                                                      \T 7.6                        
$$
    {\it  for all  \f \xi  \A   and \f z \A  satisfying  
     \f \langle  \xi  \rangle  \,  \le \,  K \A  and   \f |z|  \ge 2K$.}
\vskip 6pt
    {\it Proof.}   We prove the lemma by induction on the length
    of \f \alpha$. 
    Let \f K > 1$.  If 
    \f \langle  \xi  \rangle  \,  \le \,  K \A  and   \f |z|  \ge 2K$,
    then we see that
$$
   | \, {\langle  \xi  \rangle }^2 \, - \, z^2  \, |
       \,  \ge \,  
             \frac34 \, |z|^2                                   \T 7.7
$$
    and
$$
      | \, \widehat L_0 ( \xi ) \, + \, zI \, |
                  \,  \le \, C \, |z|.
$$
    Hence
$$
   | \, R(\xi; \, z) \, | \;   \le  \,   
          \frac{\,C\,}{|z|}  
             \qq ( \,  \langle  \xi  \rangle 
               \,  \le \,  K, \;\;   |z|  \ge 2K   \, ),        \T 7.8
$$
    which proves (7.6) for \f  \alpha = 0$.    We next prove (7.6)
    for \f \alpha \A with \f | \alpha | = 1$.
    Since
$$
\frac{\pa R}{ {\pa\xi}_j } \, = \, 
   \frac1{{\langle  \xi  \rangle }^2 \, - \, z^2}  \,
     \big( \,
      \frac{\pa \widehat L_0}{ {\pa\xi}_j } \, - \,
        2 R(\xi; \, z) \, { \xi}_j         \,     \big)  
           \qq  (j=1, \, 2, \, 3),
$$
    we deduce  from (7.7) and (7.8) that
$$
 \big| \, \frac{\pa R}{ {\pa\xi}_j } (\xi; \, z) \, \big|
    \;  \le  \;
        \frac{ C  }{\, |z|^2 \,}
     \qq ( \,  \langle  \xi  \rangle 
               \,  \le \,  K, \;\;   |z|  \ge 2K   \, ).    
                                                               \T 7.9
$$
    Here we have used the fact that \f \pa \widehat L_0 / \pa
    {\xi}_j  \A  is a constant matrix. It is
    evident that (7.9), in particular, proves (7.6) for 
    \f \alpha \A with \f | \alpha | = 1$.
                                                             \par
                                                             
     We now prove (7.6) for  \f \alpha \A with \f | \alpha | \ge 2$.
    To this end, we differentiate
    the both sides of
$$
 ( \, {\langle  \xi  \rangle }^2 \, - \, z^2 \, )  \,
    R( \xi; \, z ) \; = \;
         \widehat L_0 ( \xi ) \, + \, zI ,
                                                              \T 7.10
$$
    and apply  the Leibniz formula to the product   on the left 
    hand side of (7.10). Then we get  for 
    \f \alpha \A with \f |\alpha| \ge 2 $
$$
  \big( \frac{\pa }{\pa\xi}  \big)^{\! \alpha} \, R \; = \;
       -  \sum   \Sb  \beta \le \alpha \\
                   | \beta | = 1  \endSb   
             \binom{\alpha}{\beta} \,   2 {\xi}^{\beta} \,
     \big( \frac{\pa }{\pa\xi}  \big)^{\! \alpha - \beta}  R     \,
      - \sum   \Sb  \beta \le \alpha \\
                   | \beta | = 2  \endSb   
             \binom{\alpha}{\beta} \,   2 {\delta}_{\beta} \,
     \big( \frac{\pa }{\pa\xi}  \big)^{\! \alpha - \beta}  R,   
                                                                 \T 7.11      
$$
     where
$$
\binom{\alpha}{\beta} \,  = \,  \frac{\alpha !}
    {\beta ! ( \alpha - \beta)!}
$$
    and \f { \delta}_{\beta} = 1 \A if \f \beta \A is
    one of the following indicies \f (2, \, 0, \, 0), \;\;
    (0, \, 2, \, 0), \;\; (0, \, 0, \, 2) \A and 
    \f  {\delta}_{\beta} = 0 \A
    otherwise.  It is clear that (7.11), together with (7.8) and
    (7.9), enables us to make an 
    induction argument on \f |\alpha|$. We omit the further details.
                                                                $\square$
%

\newpage
     
        {\bf Lemma 7.3.}
\SP

\it    Let \f R( \xi; \, z) \A be the same as in Lemma 7.2. Then for
     any \f K > 1 \A and any multi-index \f \alpha \A
     there exists a constant
     \f C_{\alpha K} > 0 \A such that
$$
    \big|  \, \big( \frac{\pa}{\pa\xi} \big)^{\! \alpha} 
            \, [ \,R(\xi; z_1) \, - \, R(\xi; z_2) \, ] \, \big|
      \,  \le   \, C_{\alpha K} \, |z_1 \, - \, z_2|       
                                                                   \T 7.12 
$$
    for all \f z_1, \; z_2 \A and \f \xi \A satisfying
    \f |z_1|, \;\, |z_2|  \ge 2K \A and 
    \f \langle \xi \rangle  \le K$. 
\vskip 6pt
Proof. 
\rm     Since \f R(\xi; \, z) \, = \, ( \widehat L_0 (\xi) \, - \,
     z )^{-1}$,  we see that
$$
   R(\xi; \, z_1) \, - \, R(\xi; \, z_2) \, = \,
    (z_1 - z_2) \,   R(\xi; \, z_1) \, R(\xi; \, z_2).            \T 7.13
$$
     Applying the Leibniz formula to the right hand side of (7.13) and
     using (7.6), we get (7.12)                   
                                                                 $\square$
\vskip 12pt
     
        {\bf Lemma 7.4.}
\SP

\it    Let  \f  f \in [ \CS({\bold R}^3)]^4   \A satisfy
$$
    \text{\rm supp}[\, \widehat f \,]   \, \subset \,
\big\{  \,  \xi  \in   {\bold R}^3        \, \big/  \,\,   
                  \langle \xi \rangle 
                      \,  \le \, K   \,   \big\}                    \T 7.14
$$
    for some \f K > 1$.  For each \f z \A with \f |z| \ge 2K$, 
    define
$$
  v_z(x) \, = \, {\SF}^{\,-1} \big[ \, 
           R( \xi; \, z) \, \big] \SF  f,
$$
    where \f R( \xi; \, z) \A is the same as in Lemma 7.2.  Then
\item{\rm(i)} \it For each \f \ell \ge 0$,  there
    corresponds a constant \f C_{\ell}$, depending also on 
    \f f $, such that
$$
    |v_z|_{\ell,  \CS} \, \le \, \frac{\,C_{\ell}\,}{|z|}.
$$
\item{\rm(ii)} \it For any \f \lambda \in (-\infty, \, -2K] \cup
       [2K, \, \infty)$, 
$$
   v_{\lambda \pm i \mu}   \to  v_{\lambda}  \;\;
         \text{ in } \;\;  [ \CS({\bold R}^3)]^4   \;\;
         \text{ as } \;\;  \mu \downarrow 0.
$$

\vskip 6pt
Proof.  
\rm     Let \f \alpha \A and \f \beta \A be multi-indices. By 
    differentiation under the integral sign and 
    integration by parts, we 
    see that
$$
   x^{\alpha}\big(\frac{\pa}{\pa x}\big)^{\beta} v_z(x) 
     \, = \, (2 \pi )^{-3} 
   \int_{{\bold R}^3} \;
        e^{i x \cdot  \xi}   \;
        \big( i\,\frac{\pa}{\pa\xi}\big)^{ \! \alpha}                    
          \big\{  \, R( \xi; \, z) \, ( i \, \xi)^{\beta} \,   
                \widehat f (\xi) \,      \big\}
                                          \;   d\xi  .          \T 7.15
$$
     By Lemma 7.2, we get
$$\aligned
 {} &  | \, \text{the integrand on the RHS of (7.15)} \, |    \\
  {} & \qq \q 
     \le \; C_{\alpha \beta K} \,  \frac1{\, | z | \,} \,
        | \widehat f \,  |_{|\alpha + \beta| + 4, \CS}  \,
           \langle \xi \rangle {}^{-4}.
\endaligned
$$
    Thus we obtain
$$
  \big| \, x^{\alpha}\big(\frac{\pa}{\pa x}\big)^{\beta} v_z(x) \, \big|
     \, \le \, 
        \frac{C_{\alpha \beta K}}{|z|},
$$
     which implies conclusion (i).  Similarly, using Lemma 7.3, we can
     show that
$$
  \big| \, x^{\alpha}\big(\frac{\pa}{\pa x}\big)^{\beta}
   ( \, v_{\lambda \pm i \mu}(x) \, 
     - \, v_{\lambda}(x) \, )   \, \big|
     \, \le \, 
               C_{\alpha \beta K} \times \mu    
                                          \qq ( \mu > 0).
$$
   This leads us to conclusion (ii).                     $\square$

\vskip 12pt 
\it Proof of Theorem 2.5.   \rm   In view of Theorem 2.4 and Lemma 7.1,
    it is sufficient to show that (2.4) is true for any \f f \in
    { \CX}_0 $. Let \f f \A be in  \f { \CX}_0 \A 
    and choose \f K > 1 \A  so that (7.14) is valid. 
    Define \f v_z \A in the same
    manner as in Lemma 7.4.  Recalling (6.2), we remark that
    \f R_0 (z) \, f \, = \, v_z  \A for \f z \A 
    with \f \text{Im}\; z \not= 0$.
    Then by Lemma 7.4(ii) and Theorem 2.2,
    together with Remark 2.3(ii),  we see that
$$
   R_0^{\pm} (\lambda) f \, = \, v_{\lambda} 
      \qq ( \lambda \in (-\infty, \, -2K] \cup 
                                  [2K, \, \infty) \, ).
$$
   Moreover, by Lemma 7.4(i), we have
$$
  \M \,  R_0^{\pm} (\lambda) f  \, \M_{-s}  \, \le \,
             \frac{\, C_s \,}{|\lambda|},
$$
   which trivially implies (2.4).                      $\square$


\vskip 24pt

\heading      {\bf \S8. Proof of Theorem 2.6 }   \endheading 

\SP

     Throughout this section we assume that \f Q(x) \A satisfies
     Assumption 2.1.

\vskip 12pt
     
        {\bf Lemma 8.1.}
\SP

\it    Suppose that \f 1/2 < s < (1+ \epsilon )/2$. Then there
     exists a constant \f C_* > 0 \A such that
$$
 \M \, Q \, R_0 (z) \, f \, \M_s \; 
    \le \;  C_* \,  \M \, f \, \M_s                              \T 8.1
$$
    for all  \f f  \in    \CL_{2,s} \A  and   all \f z \in J$,
    where \f J \A is the set introduced in 
    the beginning  of the proof of Theorem 6.1.
\vskip 6pt
Proof.  
\rm     Since \f s-1-\epsilon < -s$, we have
$$
 \M \, Q \,  f \, \M_s \; 
    \le \;  C_1  \,  \M \, f \, \M_{-s} 
            \qq (  f  \in    \CL_{2,-s} ) ,                    \T 8.2
$$
       where \f C_1 \A is a constant depending only on the constant
       \f K \A 
     appearing in the Assumption 2.1. Combining (8.2) with 
     Theorem 6.1 gives the lemma.                          $\square$   
\vskip 12pt
     Writing
$$
       R_t (z) \, = \, ( H_t \, - \, z)^{-1}    
$$
    we see that
$$
   R_t (z) \, \bigl(\,I \, + \, t\, Q \, R_0(z)  \bigr)  \, = \, R_0(z)
         \;\;\; \text{ on } \CL_2.
                                                               \T 8.3
$$
     According to Lemma 8.1, we can regard 
      \f I \, + \, t\, Q \,R_0(z) \A
     as a bounded operator in \f    \CL_{2,s} \A  provided that
     \f 1/2 < s < (1 + \epsilon) /2$.  
     Letting \f C_* \A be the constant  in (8.1) and choosing 
     \f t_0 > 0 \A so that
$$
     0 <  t_0 \, C_*  \; < \; 1,                                   \T 8.4
$$
    we can construct the inverse of \f I \, + \, t\,Q \,R_0(z) \A
    by using the Neumann series in  
    \f \bB(\CL_{2,s}, \CL_{2, s})  \A
    if \f |t| \le t_0 \A and \f z \in J$:
$$
  \bigl( \, I \, + \, t\,Q\,R_0(z) \, \bigr)^{-1}  \, = \,
     \sum_{\ell = 0}^{\infty}
              \bigl( \, -t \, Q \, R_0(z) \bigr)^{\ell}  .         \T 8.5    
$$
\vskip 12pt
     {\it  Proof of Theorem 2.6(i).}  \rm  
     We may assume, without loss of generality, that \f 1/2 < s 
     < (1+\epsilon)/2$.  Indeed, if 
     \f 1/2 < s < s^{\prime}$,  then
$$
 \M \, R_t^{\pm} (\lambda) \, \M_{(s^{\prime}, -s^{\prime}) }
     \, \le \,
        \M \, R_t^{\pm} (\lambda) \, \M_{(s, -s ) }.
$$

    In view of Theorem 2.2, it is sufficient to show that
    there exists a constant \f C_2 > 0 \A  such that
$$
     \M \,  R_t (z) \, \M_{(s, -s)} \;  \le  \;    
                                C_2                           \T 8.6
$$
     for all \f t \A  with 
      \f |t| \le t_0 \A and all \f z \in J$.
           By (8.5) and Lemma 8.1, we see that
$$
  \big\M \,  \bigl( \, I \, + \, t \, Q \, R_0(z) \,  \bigr)^{-1}  \,
       \big\M_{(s,s)}   \,  \le   \, 
          \frac1{\, 1 \, - \, t_0 C_* \,} 
                 \qq ( |t| \le t_0, \; \;  z \in J ).            \T 8.7
$$
    Hence, using (8.3) and (8.7), we have
$$
   \M \, R_t (z) \, \M_{(s, -s)}  \,  \le \,
              \M \, R_0 (z) \, \M_{(s, -s)}         \, 
                         \frac1{\, 1 \, - \, t_0 C_* \,} 
$$
    when \f |t|  \le t_0 \A and \f z \in J$.  Combining  this
    inequality with Theorem 6.1, we obtain  (8.6).             $\square$ 

\vskip 12pt
     {\it  Proof of Theorem 2.6(ii).}  \rm  
     We may assume again, without loss of generality, that \f 1/2 < s 
     < (1+\epsilon)/2$.  In fact, if 
     \f 1/2 < s < s^{\prime}$,  then
     \f  \M \, f  \, \M_{-s^{\prime}}
     \, \le \,
        \M \, f \, \M_{ -s}$.

        Suppose that \f   f  \in    \CL_{2,s}$,  \f z \in J$   and
    \f |t| \le t_0$.  Then, by (8.3) and (8.5), we have
$$\aligned
 \M \, R_t (z) \, f \, \M_{-s}  \, & \le \,
        \M \, R_0 (z) \, f \,  \M_{-s}                           
       \;  + \; 
       \M \, R_0 (z) \, \sum_{\ell = 1}^N 
              \bigl( \, -t \, Q \, R_0(z)  \bigr)^{\ell}  
                     \,    f \,  \M_{-s}                           \\
       {} & \q \;\;  + \; 
       \M \, R_0 (z) \, \sum_{\ell = N+1}^{\infty} 
              \bigl( \, -t \, Q \, R_0(z)  \bigr)^{\ell}  
                     \,    f \,  \M_{-s}                      
\endaligned                                                       \T 8.8
$$
    for any positive integer \f N$. The second term on the right
    hand side of (8.8)
    is estimated by
$$\split
  \M \, R_0 (z) \, & \M_{(s, -s)} \,  \
        \big\{ 
      \sum_{\ell = 1}^N 
          \big\M \,  \bigl( \, -t \, Q \, R_0(z)  \bigr)^{\ell - 1}   \,
                   \big\M_{(s,s)} \,    \big\}
                     \,  |t|
                      \;  \M \,  Q \, R_0(z)   f \,  \M_s           \\
       {} &  \le \;    C_3  \,
           \big\{
         \sum_{\ell = 1}^N    ( t_0 \, C_* )^{\ell -1}    \big\}
           \, t_0 \, C_1  \,  \M \,  R_0(z)   f \,  \M_{-s}
\endsplit
$$
    where \f C_1 \A is the same constant as in (8.2)  and \f C_3   \A
    is the constant given by the supremum  in Theorem 6.1.
    With this notation, the last term on the right hand side of (8.8)
    is less than or equal to
$$
 C_3  \,
           \big\{
         \sum_{\ell = N + 1}^{\infty} ( t_0 \, C_* )^{\ell}    \big\}
       \,  \M \,  f \,  \M_s.
$$
    Summing up, we get
$$\aligned
 \M \, R_t (z) \, f \, \M_{-s}  \, & \le \,
    \big\{  \, 1 \, +  \, t_0 \, C_1 \, C_3 
       \sum_{\ell = 1}^N ( t_0 \, C_* )^{\ell-1}    \big\}
        \,  \M \, R_0 (z) \, f \,  \M_{-s}                         \\
   {} & \qq  \q  + \;   
      C_3 \,   \big\{ 
       \sum_{\ell = N + 1}^{\infty} ( t_0 \, C_* )^{\ell}  \big\}
        \,  \M \, f \,  \M_s       
\endaligned                                                       \T 8.9
$$
    for any positive integer \f N$. 
    We now replace  \f z \A in (8.9) with
      \f  \lambda \, \pm \, i \, \mu \A  ( \f |\lambda| > 2$, 
        \f 0 < \mu <1 \, $) and take the limits as 
        \f \mu \downarrow 0$. Then we obtain, by Theorem 2.2,
$$\aligned
 \M \, R_t^{\pm} (\lambda) \, f \, \M_{-s}  \, & \le \,
    \big\{  \, 1 \, +  \, t_0 \, C_1 \, C_3 
       \sum_{\ell = 1}^N ( t_0 \, C_* )^{\ell-1}    \big\}
        \,  \M \, R_0^{\pm} (\lambda) \, f \,  \M_{-s}                         \\
   {} & \qq  \q + \;  
      C_3 \,   \big\{ 
       \sum_{\ell = N + 1}^{\infty} ( t_0 \, C_* )^{\ell}  \big\}
        \,  \M \, f \,  \M_s       
\endaligned                                                       \T 8.10
$$
    for any positive integer \f N$.   Theorem 2.5, together with (8.10),
    implies that 
$$
\limsup_{|\lambda| \to \infty} \,
     \M \, R_t^{\pm} (\lambda) \, f \, \M_{-s}  \,  \le \,
  C_3 \,   \big\{ 
       \sum_{\ell = N + 1}^{\infty} ( t_0 \, C_* )^{\ell}  \big\}
        \,  \M \, f \,  \M_s. 
$$
    Since \f N \A is arbitrary and  \f 0 < t_0 \, C_* <1 \A 
    ( recall (8.4) ), we conclude that
$$
\lim_{|\lambda| \to \infty} \,
     \M \, R_t^{\pm} (\lambda) \, f \, \M_{-s}  \,  = \, 0
$$
    for  \f   f  \in    \CL_{2,s}$  and  \f t \, \A with
    \f |t| \le t_0$.                                          $\square$


\newpage   

\centerline{ {\bf References} }
\BP

\item{[1]} \ S. Agmon,  Spectral properties of 
     Schr\"odinger operators and scattering theory,
     Ann. Scoula Norm. Sup. Pisa ({\bf 4})2 (1975), 151-218. 
\item{[2]} \ A.P. Calder\'on and R. Vaillancourt,  A class of
     bounded pseudo-differential operators, Proc. Nat. Acad. Sci.
     USA {\bf 69} (1972), 1185-1187.
\item{[3]} \ T. Ikebe and Y. Sait$\bar{ \hbox{o}}$, Limiting
     absorption method and absolute continuity for the 
     Schr\"odinger operators, \rm  J. Math. Kyoto
     Univ. {\bf 7} (1972), 513-542.
\item{[4]} \ T. Kato, {\it  Perturbation Theory for Linear Operators, 2nd 
     ed.,} Springer-Verlag, 1976.
\item{[5]} \ H. Kumano-go, {\it Pseudo-differential Operators }, the MIT
     Press, 1981.
\item{[6]} \ C. Pladdy, Y. Sait\B o and T. Umeda, 
     Asymptotic behavior of the resolvent of the Dirac operator,
     to appear in the Proceedings of International Conference on
     Mathematical Results in Quantum Mechanics in Blossin, Germany
     (May 17--23, 1993) to be published in the series \lq\lq Operator
     Theory: Advances and Applications" (ed. by I. Gohberg).
\item{[7]} \ C. Pladdy, Y. Sait\B o and T. Umeda, Radiation 
     condition for Dirac operators, preprint.     
\item{[8]} \  Y. Sait$\bar{ \hbox{o}}$, The principle of
     limiting absorption for second-order differential operators
     with operator-valued coefficients, Publ. Res. Inst. Math. Sci. Kyoto 
     Univ. {\bf 7} (1972), 581-619.
\item{[9]} \ Y. Sait$\bar{ \hbox{o}}$,  The principle of limiting absorption 
     for the non-selfadjoint Schr\"odin\-ger
       operator in ${\bold R}^N  \,\,\, 
     (N\not= 2)$, Publ. Res. Inst. Math. Sci. Kyoto Univ. {\bf 9} (1974),
     397-428.
\item{[10]} \ Y. Sait$\bar{ \hbox{o}}$, The principle of limiting absorption 
     for the non-selfadjoint Schr\"odin\-ger operator in ${\bold R}^2$,
      Osaka J.
     Math.  {\bf 11} (1974), 295-306.
\item{[11]} \ Y. Sait$\bar{ \hbox{o}}$, An asymptotic behavior of
      the $\CS$--matrix and the inverse scattering problem, J. Math.
      Phys. {\bf 25} (1984), 3105-3111.
\item{[12]} \ Y. Sait$\bar{ \hbox{o}}$, An  approximation formula
      in the inverse scattering problem, J. Math. Phys. {\bf 27} 
      (1986), 1145-1153.      
\item{[13]} \ O. Yamada, On the principle of limiting absorption for the 
     Dirac operators, Publ. Res. Inst. Math. Sci. Kyoto Univ. {\bf 8} 
     (1972/73), 557-577.
\item{[14]} \ O. Yamada, Eigenfunction expansions and scattering
     theory for Direc operators,  Publ. Res. Inst. Math. Sci. Kyoto Univ.
      {\bf 11} (1976), 651-689. 
\item{[15]} \ O. Yamada, A remark on the limiting absorption method
     for Dirac operators, preprint, 1993.
   
\enddocument